\documentclass[11pt, english]{article}

\usepackage[english]{babel}
\usepackage{mypreprint}
\usepackage[font=small]{subfig} 
\usepackage{units} 
\usepackage{graphicx}
\usepackage{epstopdf}
\usepackage{amsmath,amssymb}
\usepackage{xcolor}
\usepackage[numbers]{natbib} 
\usepackage{tikz,pgfplots}
\pgfplotsset{compat=newest}
\usepackage{pgfkeys}
\newlength\myheight
\newlength\mywidth

\usetikzlibrary{external} 
\usepgfplotslibrary{external} 
\definecolor{tumblue}{HTML}{0065BD}
\definecolor{tumblue1}{HTML}{98C6EA}
\definecolor{tumblue2}{HTML}{64A0C8}
\definecolor{tumblue3}{HTML}{0073CF}
\definecolor{tumblue4}{HTML}{005293}
\definecolor{tumblue5}{HTML}{003359}
\definecolor{tumgreen}{HTML}{A2AD00}
\definecolor{tumorange}{HTML}{E37222}
\definecolor{tumivory}{HTML}{DAD7CB}
\definecolor{tumviolet}{HTML}{69085A}
\definecolor{tumred}{HTML}{C4071B}
\definecolor{tumgray}{rgb}{0.4,0.4,0.4}
\definecolor{tumgreen2}{rgb}{.15,.50,.20} 
\definecolor{tumbluegreen}{rgb}{0,.37,.44} 

\usepackage{hyperref}

\begin{document}

\date{June 15, 2016}

\title{Model reduction of linear time-varying systems with applications for moving loads}

\author{Maria Cruz Varona\footnotemark[2] \
and Boris Lohmann\footnotemark[2]}

\maketitle

\renewcommand{\thefootnote}{\fnsymbol{footnote}}
\footnotetext[2]{Chair of Automatic Control,
				Technical University of Munich,
				Boltzmannstr. 15, D-85748 Garching (\textbraceleft \texttt{maria.cruz,lohmann}\textbraceright\texttt{@tum.de})}
\renewcommand{\thefootnote}{\arabic{footnote}}

\begin{abstract}
In this paper we consider different model reduction techniques for systems with moving loads. Due to the time-dependency of the input and output matrices, the application of time-varying projection matrices for the reduction offers new degrees of freedom, which also come along with  some challenges. This paper deals with both simple methods for the reduction of particular linear time-varying systems, as well as with a more advanced technique considering the emerging time derivatives.

\end{abstract}

\section{Introduction} \label{sec:1}
\noindent The detailed modeling of physical and technical phenomena arising in many engineering and computer science applications may yield models of very large dimension. This is particular the case in fields such as thermo-fluid dynamics, structural mechanics or integrated circuit design, where the models are mostly obtained from a spatial discretization of the underlying partial differential equations. The resulting large systems of ordinary differential equations or differential-algebraic equations are computationally expensive to simulate and handle. In order to reduce the computational effort, model reduction techniques that generate reduced-order models that approximate the dynamic behaviour and preserve the relevant properties of the original model are required. For the reduction of linear time-invariant (LTI) systems, various well-established reduction approaches exist (see e.g. \cite{Cru_Antoulas_Book}). In the past ten years, further model reduction methods have been developed for linear, parametric and nonlinear systems \cite{Cru_Benner_2015_Survey,Cru_Baur_2014_Survey} and applied in a wide variety of domains.

In this contribution, we investigate model order reduction of linear time-varying (LTV) systems. Such systems arise in many real-life applications, since dynamical systems often depend on parameters which vary over time or might alter their behaviour due to ageing, degradation, environmental changes and time-dependent operating conditions. Another possible application for LTV systems are moving loads. This particular but still very frequent problem arises, for example, in working gears, cableways, bridges with moving vehicles or milling processes. Since the position of the acting force varies over time, systems with sliding components exhibit a time-variant behaviour. The varying load location can be modeled and considered in different ways, thus yielding diverse alternative representations for systems with moving loads and, according to this, leading to different approaches to reduce them.

One possibility is to represent moving loads as LTV systems, in which \emph{only} the input and/or output matrices are time-dependent. Such systems can be then reduced using balanced truncation model reduction methods developed in \cite{Cru_Shokoohi_1983_LTV-BT, Cru_Sandberg_2004_LTV-BT}. These approaches, however, require a high computational and storage effort, since two differential Lyapunov equations must be solved. Recently, a practical and efficient procedure of balanced truncation for LTV systems has been presented in \cite{Cru_Lang_MCMDS_2016}. Note that these aforementioned balanced truncation techniques can be applied to general LTV systems, where all system matrices are time-dependent. For the reduction of systems with only time-varying input and output matrices the two-step approach proposed in \cite{Cru_Stykel_CAM_2016,Cru_Baumann_Comparison_2016} can also be pursued. This method consists first on a low-rank approximation of the time-dependent input matrix and consequently on applying standard model reduction techniques to the resulting LTI system with a modified input. The approximation of the input matrix in a low-dimensional subspace is performed via the solution of a linear least squares minimization problem. 

Systems with moving loads can further be modeled by means of linear switched systems. Well-known reduction methods such as balanced truncation can then be applied for the reduction of each LTI subsystem \cite{Cru_Lang_2014}.

A last alternative option for describing systems with moving loads is to consider the load position as a time-dependent parameter of the system model. This results in a linear parameter-varying (LPV) system, in which only the input and/or output matrices depend on a time-varying parameter. In many recent publications, e.g. \cite{Cru_Fischer_2014_Application,Cru_Lang_2014,Cru_Fischer_ECCOMAS_2015,Cru_Baumann_Comparison_2016}, the parameter is assumed to be time-independent. Thereby any parametric model order reduction (pMOR) approach \cite{Cru_Benner_2015_Survey} can be applied to the resulting parametric LTI system. In some other recent publications \cite{Cru_Tamarozzi_2014,Cru_Cruz_MATHMOD_2015,Cru_Cruz_Hirschberg_2015,Cru_Cruz_MoRePaS_2015} the time variation of the parameter is taken into account, whereby new time derivative terms emerge during the time-dependent parametric model reduction process.

This paper deals with different time-varying model reduction techniques for systems with moving loads. Firstly, LTV systems are considered and the time-dependent projective reduction framework is briefly explained in section~\ref{sec:2}. Since moving loads represent particular LTV systems, we then introduce some straightforward reduction approaches for the resulting special cases in section~\ref{sec:3}. In the second part of the paper, we focus on LPV systems and present a time-dependent parametric model reduction approach by matrix interpolation \cite{Cru_Cruz_MATHMOD_2015,Cru_Cruz_MoRePaS_2015} in section~\ref{sec:4}. Some numerical results for the reduction of systems with moving loads applying the proposed methods are finally reported and discussed in section~\ref{sec:5}.

\section{Linear Time-Varying Model Order Reduction} \label{sec:2}
\noindent In the following we first consider a high-dimensional linear time-varying system of the form

\begin{equation}
	\begin{aligned}
		\mathbf{E}(t) \, \dot{\mathbf{x}}(t) &= \mathbf{A}(t) \, \mathbf{x}(t) + \mathbf{B}(t) \, \mathbf{u}(t),\\[0.2em]
		\mathbf{y}(t) &= \mathbf{C}(t) \, \mathbf{x}(t), 
	\end{aligned}
	\label{equ:LTV_system}
\end{equation}
where $\mathbf{E}(t)$, $\mathbf{A}(t) \in \mathbb{R}^{n \times n}$, $\mathbf{B}(t) \in \mathbb{R}^{n \times m}$ and $\mathbf{C}(t) \in \mathbb{R}^{q \times n}$ are the time-dependent system matrices, $\mathbf{x}(t) \in \mathbb{R}^n$ is the state vector and $\mathbf{u}(t) \in \mathbb{R}^m$, $\mathbf{y}(t) \in \mathbb{R}^q$ represent the inputs and outputs of the system, respectively. The system matrix $\mathbf{E}(t)$ is assumed to be nonsingular for all $t \in \left[0, T\right]$. Note that it is also possible to consider second-order systems and reformulate them into the first-order form \eqref{equ:LTV_system}.

\subsection{Time-dependent Projective Reduction Framework} \label{subsec:tMOR}
\noindent In projective model order reduction, we aim to find a reduced-order model by approximating the state vector $\mathbf{x}(t)$ on a subspace of lower dimension $r \ll n$. In the time-varying case, the state vector $\mathbf{x}(t)$ might be projected onto a \emph{varying} subspace spanned by the columns of a \emph{time-dependent} projection matrix $\mathbf{V}(t) \in \mathbb{R}^{n \times r}$ \cite{Cru_Shokoohi_1983_LTV-BT,Cru_Tamarozzi_2014}. Therefore, the approximation equations read
\begin{equation}
	\begin{aligned}
		\mathbf{x}(t) &\approx \mathbf{V}(t) \, \mathbf{x}_r(t), \\[0.2em]
		\dot{\mathbf{x}}(t) &\approx \dot{\mathbf{V}}(t) \, \mathbf{x}_r(t) + \mathbf{V}(t) \, \dot{\mathbf{x}}_r(t),
	\end{aligned}
\end{equation}
whereby the product rule must be considered in this case for the time derivative of the state vector. Plugging first these both equations into \eqref{equ:LTV_system}, and applying thereon a properly chosen time-dependent projection matrix $\mathbf{W}(t)$ which enforces the Petrov-Galerkin condition leads to the time-varying reduced-order model
\begin{equation}
	\resizebox{0.93\linewidth}{!}{
		$\begin{aligned}
		\overbrace{\mathbf{W}(t)^T \mathbf{E}(t) \mathbf{V}(t)}^{\mathbf{E}_r(t)} \, \dot{\mathbf{x}}_r(t) &= 
		\left(\overbrace{\mathbf{W}(t)^T \mathbf{A}(t) \mathbf{V}(t)}^{\mathbf{A}_r(t)} - \mathbf{W}(t)^T \mathbf{E}(t) 
		\textcolor{tumred}{\dot{\mathbf{V}}(t)}  \right) \, \mathbf{x}_r(t) + 
		\overbrace{\mathbf{W}(t)^T \mathbf{B}(t)}^{\mathbf{B}_r(t)} \, \mathbf{u}(t), \\[0.2em]
		\mathbf{y}_r(t) &= \underbrace{\mathbf{C}(t) \mathbf{V}(t)}_{\mathbf{C}_r(t)} \, \mathbf{x}_r(t).
		\end{aligned}$
	}
	\label{equ:ROM_LTV_system}
\end{equation}
It is noteworthy to mention that the system matrix of the reduced-order model \eqref{equ:ROM_LTV_system} not only comprises the reduced matrix $\mathbf{A}_r(t)$, but also includes a further term which depends on the time derivative $\dot{\mathbf{V}}(t)$ of the time-varying projection matrix. This additional term influences the dynamic behaviour of the reduced-order model and should therefore be taken into account.

The usage of time-dependent projection matrices for the reduction of linear time-varying systems certainly offers benefits regarding the approximation quality. For their computation, however, standard reduction methods such as balanced truncation cannot be directly applied, but must be adapted instead. Furthermore, the time derivative of $\mathbf{V}(t)$ should be approximated numerically (thus increasing the computational effort) and included in the time integration scheme of the reduced-order model \cite{Cru_Lang_MCMDS_2016}.

\section{Straightforward Reduction Approaches for particular Linear Time-Varying Systems} \label{sec:3}
\noindent In the previous section we have seen that the application of time-dependent projection matrices for the reduction of LTV systems comes along with some difficulties and challenges. For the reduction of particular LTV systems, in which only the input and/or output matrices depend on time, the usage of \emph{time-independent} projection matrices $\mathbf{V}$ and $\mathbf{W}$ might be sufficient. In this section we discuss some special cases for LTV systems and propose straightforward approaches to reduce them. 
\subsection{Case 1: Moving Loads}
\noindent The first case we want to consider is a high-dimensional LTV system with only time-varying input matrix, and all other matrices being time-independent:
\begin{equation}
	\begin{aligned}
		\mathbf{E} \, \dot{\mathbf{x}}(t) &= \mathbf{A} \, \mathbf{x}(t) + \textcolor{tumorange}{\mathbf{B}(t)} \, \mathbf{u}(t),\\
		\mathbf{y}(t) &= \mathbf{C} \, \mathbf{x}(t). 
	\end{aligned}
	\label{equ:LTV_B(t)}
\end{equation}
The time-dependent input matrix describes the position of the moving forces at time $t$. In the following we present two straightforward approaches to reduce a system in the form above using time-independent projection matrices.
\subsubsection*{Approach 1: Two-step method}
\noindent The first straightforward reduction method is deducted from the two-step approach presented in \cite{Cru_Stykel_CAM_2016,Cru_Baumann_Comparison_2016}. The method is composed of two steps:

\begin{enumerate}
	\item The time-variability of the input matrix is shifted to the input variables through a low-rank approximation of the input matrix by $\textcolor{tumorange}{\mathbf{B}(t)} \approx \textcolor{tumviolet}{\mathbf{B}} \ \tilde{\mathbf{B}}(t)$, where $\textcolor{tumviolet}{\mathbf{B}} \in \mathbb{R}^{n \times \tilde{m}}$ with $\tilde{m} \ll n$ is a constant matrix and $\tilde{\mathbf{B}}(t) \in \mathbb{R}^{\tilde{m} \times m}$.
	Introducing a new input $\tilde{\mathbf{u}}(t)=\tilde{\mathbf{B}}(t) \, \mathbf{u}(t)$, the original model \eqref{equ:LTV_B(t)} can be transformed to:
	\begin{equation}
		\begin{aligned}
			\mathbf{E} \, \dot{\mathbf{x}}(t) &= \mathbf{A} \, \mathbf{x}(t) + \textcolor{tumviolet}{\mathbf{B}} \, \overbrace{\tilde{\mathbf{B}}(t) \, \mathbf{u}(t)}^{\tilde{\mathbf{u}}(t)},\\
			\mathbf{y}(t) &= \mathbf{C} \, \mathbf{x}(t).
		\end{aligned}
	\end{equation}
	\item The resulting multiple-input multiple-output (MIMO) LTI system $(\mathbf{E}, \mathbf{A}, \textcolor{tumviolet}{\mathbf{B}}, \mathbf{C})$ can subsequently be reduced by means of balanced truncation, MIMO rational Krylov or MIMO-IRKA, for instance. The reduced-order model is then given by
	\begin{equation}
		\begin{aligned}
			\overbrace{\mathbf{W}^T \mathbf{E} \mathbf{V}}^{\mathbf{E}_r} \, \dot{\mathbf{x}}_r(t) &= \overbrace{\mathbf{W}^T \mathbf{A} \mathbf{V}}^{\mathbf{A}_r} \, \mathbf{x}_r(t) + \overbrace{\mathbf{W}^T \textcolor{tumviolet}{\mathbf{B}}}^{\mathbf{B}_r} \, \overbrace{\tilde{\mathbf{B}}(t) \, \mathbf{u}(t)}^{\tilde{\mathbf{u}}(t)},\\
			\mathbf{y}_r(t) &= \underbrace{\mathbf{C} \mathbf{V}}_{\mathbf{C}_r} \, \mathbf{x}_r(t), 
		\end{aligned}
	\end{equation}	
	where the reduced time-varying input matrix reads $\mathbf{B}_r(t)=\mathbf{B}_r \, \tilde{\mathbf{B}}(t)$.
\end{enumerate}
For the approximation of the input matrix $\mathbf{B}(t)$, other than \cite{Cru_Stykel_CAM_2016,Cru_Baumann_Comparison_2016} we simply take the correct input columns $\mathbf{b}_i(t)$ with the moving load acting at corresponding nodes $i$ of a coarse finite element grid and form the low-rank matrix $\mathbf{B}$ with them, without performing a least squares minimization with the basis functions. Note that the two-step approach only provides satisfactory results, if the number of columns $\tilde{m}$ of the low-rank matrix $\mathbf{B}$ is sufficiently large \cite{Cru_Stykel_CAM_2016}. Otherwise the overall approximation error in the output (due to the approximation error in the input matrix and the model reduction error) can become inadmissibly large. Note also that this reduction method is limited to systems with a known trajectory of the load before the simulation.
\subsubsection*{Approach 2: One-sided reduction with output Krylov subspace}
\noindent The second straightforward method uses Krylov subspaces for the reduction and exploits the fact that the only time-varying element in system \eqref{equ:LTV_B(t)} is the input matrix $\mathbf{B}(t)$. Since an input Krylov subspace would yield a time-varying projection matrix
\begin{equation}
	\begin{aligned}
		\mathbf{V}(t) := \begin{bmatrix}
			{\mathbf{A}}^{-1}_{s_0}\textcolor{tumorange}{\mathbf{B}(t)} & {\mathbf{A}}^{-1}_{s_0}\mathbf{E}{\mathbf{A}}^{-1}_{s_0}\textcolor{tumorange}{\mathbf{B}(t)} & \ldots & ({\mathbf{A}}^{-1}_{s_0}\mathbf{E})^{r-1}{\mathbf{A}}^{-1}_{s_0}\textcolor{tumorange}{\mathbf{B}(t)} 
		\end{bmatrix},
	\end{aligned}
	\label{equ:V(t)_B(t)}
\end{equation}
where $\mathbf{A}_{s_0}=\mathbf{A} - s_0 \mathbf{E}$, the idea of this approach is to perform a one-sided reduction with $\mathbf{V} = \mathbf{W}$, where the columns of $\mathbf{W}$ form a basis of the output Krylov subspace:
\begin{equation}
	\begin{aligned}
		\mathbf{W} := \begin{bmatrix}
			{\mathbf{A}}^{-T}_{s_0} \mathbf{C}^T & {\mathbf{A}}^{-T}_{s_0}\mathbf{E}^T{\mathbf{A}}^{-T}_{s_0}\mathbf{C}^T & \ldots & ({\mathbf{A}}^{-T}_{s_0}\mathbf{E}^T)^{r-1}{\mathbf{A}}^{-T}_{s_0} \mathbf{C}^T 
		\end{bmatrix}.
	\end{aligned}
\end{equation}
Thereby, time-independent projection matrices are obtained for computing the reduced-order model
\begin{equation}
	\begin{aligned}
		\overbrace{\mathbf{W}^T \mathbf{E} \mathbf{W}}^{\mathbf{E}_r} \, \dot{\mathbf{x}}_r(t) &= \overbrace{\mathbf{W}^T \mathbf{A} \mathbf{W}}^{\mathbf{A}_r} \, \mathbf{x}_r(t) + \overbrace{\mathbf{W}^T \textcolor{tumorange}{\mathbf{B}(t)}}^{\mathbf{B}_r(t)} \, \mathbf{u}(t),\\
		\mathbf{y}_r(t) &= \underbrace{\mathbf{C} \mathbf{W}}_{\mathbf{C}_r} \, \mathbf{x}_r(t). 
	\end{aligned}
\end{equation}
Although only the first $r$ Taylor coefficients (so-called \emph{moments}) of the transfer function of the original and the reduced model around the expansion points $s_0$ match due to the application of a one-sided reduction, we obtain time-independent projection matrices with this approach and can therefore get rid of the time derivative $\dot{\mathbf{V}}(t)$.
\subsection{Case 2: Moving Sensors}
\noindent Now we consider a LTV system with only time-varying output matrix

\begin{equation}
	\begin{aligned}
		\mathbf{E} \, \dot{\mathbf{x}}(t) &= \mathbf{A} \, \mathbf{x}(t) + \mathbf{B} \, \mathbf{u}(t),\\
		\mathbf{y}(t) &= \textcolor{tumblue}{\mathbf{C}(t)} \, \mathbf{x}(t). 
	\end{aligned}
	\label{equ:LTV_C(t)}
\end{equation}
The time-dependent output matrix describes the position of the moving sensors at time $t$. This particular LTV system can easily be reduced in the following ways.  
\subsubsection*{Approach 1: Two-step method}
\begin{enumerate}
	\item We shift the time-variability of the output matrix to the output variables through a low-rank approximation by $\textcolor{tumblue}{\mathbf{C}(t)} \approx \tilde{\mathbf{C}}(t) \, \textcolor{tumviolet}{\mathbf{C}}$, where $\textcolor{tumviolet}{\mathbf{C}} \in \mathbb{R}^{\tilde{q} \times n}$ with $\tilde{q} \ll n$ is a constant matrix and $\tilde{\mathbf{C}}(t) \in \mathbb{R}^{q \times \tilde{q}}$.
	Introducing a new output $\tilde{\mathbf{y}}(t)=\mathbf{C} \, \mathbf{x}(t)$, the original model \eqref{equ:LTV_C(t)} can be transformed to:
	\begin{equation}
		\begin{aligned}
			\mathbf{E} \, \dot{\mathbf{x}}(t) &= \mathbf{A} \, \mathbf{x}(t) + \mathbf{B} \, \mathbf{u}(t),\\
			\mathbf{y}(t) &= \tilde{\mathbf{C}}(t) \, \underbrace{\textcolor{tumviolet}{\mathbf{C}} \, \mathbf{x}(t)}_{\tilde{\mathbf{y}}(t)}.
		\end{aligned}
	\end{equation}
	\item The resulting system $(\mathbf{E}, \mathbf{A}, \mathbf{B}, \textcolor{tumviolet}{\mathbf{C}})$ can subsequently be reduced by means of any appropriate multiple-input multiple-output LTI reduction technique. The calculated time-independent projection matrices lead to the reduced-order model
	\begin{equation}
		\begin{aligned}
			\overbrace{\mathbf{W}^T \mathbf{E} \mathbf{V}}^{\mathbf{E}_r} \, \dot{\mathbf{x}}_r(t) &= \overbrace{\mathbf{W}^T \mathbf{A} \mathbf{V}}^{\mathbf{A}_r} \, \mathbf{x}_r(t) + \overbrace{\mathbf{W}^T \mathbf{B}}^{\mathbf{B}_r} \, \mathbf{u}(t),\\
			\mathbf{y}_r(t) &= \tilde{\mathbf{C}}(t) \, \underbrace{\textcolor{tumviolet}{\mathbf{C}} \mathbf{V}}_{\mathbf{C}_r} \, \mathbf{x}_r(t),
		\end{aligned}
	\end{equation}
	with the reduced time-varying output matrix $\mathbf{C}_r(t) = \tilde{\mathbf{C}}(t) \, \mathbf{C}_r$.
\end{enumerate}
The approximation of $\mathbf{C}(t)$ is performed by simply taking the output rows with the moving sensor at the corresponding nodes of a coarse finite element grid. Note that the approximation of the output matrix yields additional errors in the output \cite{Cru_Stykel_CAM_2016}.
\subsubsection*{Approach 2: One-sided reduction with input Krylov subspace}
\noindent Since in this case an output Krylov subspace would lead to a time-varying projection matrix 
\begin{equation}
	\begin{aligned}
		\mathbf{W}(t) := \begin{bmatrix}
			{\mathbf{A}}^{-T}_{s_0}\textcolor{tumblue}{\mathbf{C}(t)^T} & {\mathbf{A}}^{-T}_{s_0}\mathbf{E}^T{\mathbf{A}}^{-T}_{s_0}\textcolor{tumblue}{\mathbf{C}(t)^T} & \ldots & ({\mathbf{A}}^{-T}_{s_0}\mathbf{E}^T)^{r-1}{\mathbf{A}}^{-T}_{s_0}\textcolor{tumblue}{\mathbf{C}(t)^T} 
		\end{bmatrix}
	\end{aligned}
	\label{equ:W(t)_C(t)}
\end{equation}
due to the time-dependent output matrix $\mathbf{C}(t)$, the idea is now to perform a one-sided reduction with $\mathbf{W} = \mathbf{V}$, where the columns of $\mathbf{V}$ form a basis of the input Krylov subspace:
\begin{equation}
	\begin{aligned}
		\mathbf{V} := \begin{bmatrix}
			{\mathbf{A}}^{-1}_{s_0}\mathbf{B} & {\mathbf{A}}^{-1}_{s_0}\mathbf{E}{\mathbf{A}}^{-1}_{s_0}\mathbf{B} & \ldots & ({\mathbf{A}}^{-1}_{s_0}\mathbf{E})^{r-1}{\mathbf{A}}^{-1}_{s_0}\mathbf{B}
		\end{bmatrix}.
	\end{aligned}
\end{equation}
The reduced model is then given by:
\begin{equation}
	\begin{aligned}
		\overbrace{\mathbf{V}^T \mathbf{E} \mathbf{V}}^{\mathbf{E}_r} \, \dot{\mathbf{x}}_r(t) &= \overbrace{\mathbf{V}^T \mathbf{A} \mathbf{V}}^{\mathbf{A}_r} \, \mathbf{x}_r(t) + \overbrace{\mathbf{V}^T \mathbf{B}}^{\mathbf{B}_r} \, \mathbf{u}(t),\\
		\mathbf{y}_r(t) &= \underbrace{\textcolor{tumblue}{\mathbf{C}(t)} \mathbf{V}}_{\mathbf{C}_r(t)} \, \mathbf{x}_r(t). 
	\end{aligned}
\end{equation}
Due to the application of a one-sided reduction, only $r$ moments are matched. Nevertheless, the time derivative is avoided, since $\mathbf{V}$ and $\mathbf{W}$ are time-independent ($\dot{\mathbf{V}} = \mathbf{0}$).
\subsection{Case 3: Moving Loads and Sensors}
\noindent Finally, we consider the combined case with time-varying input \emph{and} output matrices
\begin{equation}
	\begin{aligned}
		\mathbf{E} \, \dot{\mathbf{x}}(t) &= \mathbf{A} \, \mathbf{x}(t) + \textcolor{tumorange}{\mathbf{B}(t)} \, \mathbf{u}(t),\\
		\mathbf{y}(t) &= \textcolor{tumblue}{\mathbf{C}(t)} \, \mathbf{x}(t). 
	\end{aligned}
\end{equation}
If the sensor position coincides with the location of the load, then $\mathbf{C}(t) = \mathbf{B}(t)^T$.
\subsubsection*{Approach 1: Two-step method}
\noindent In this case, the respective two-step techniques explained before have to be combined properly:

\begin{enumerate}
	\item The time-variability of $\mathbf{B}(t)$ is shifted to the input variables and the time-dependency of $\mathbf{C}(t)$ to the output variables, thus obtaining a MIMO LTI system.
	\item Time-independent projection matrices are then calculated applying an appropriate model order reduction method to the resulting system $(\mathbf{E}, \mathbf{A}, \textcolor{tumviolet}{\mathbf{B}}, \textcolor{tumviolet}{\mathbf{C}})$. The reduced-order model is finally given by
	\begin{equation}
		\begin{aligned}
			\overbrace{\mathbf{W}^T \mathbf{E} \mathbf{V}}^{\mathbf{E}_r} \, \dot{\mathbf{x}}_r(t) &= \overbrace{\mathbf{W}^T \mathbf{A} \mathbf{V}}^{\mathbf{A}_r} \, \mathbf{x}_r(t) + \overbrace{\mathbf{W}^T \textcolor{tumviolet}{\mathbf{B}}}^{\mathbf{B}_r} \, \overbrace{\tilde{\mathbf{B}}(t) \, \mathbf{u}(t)}^{\tilde{\mathbf{u}}(t)},\\
			\mathbf{y}_r(t) &= \tilde{\mathbf{C}}(t) \, \underbrace{\textcolor{tumviolet}{\mathbf{C}} \mathbf{V}}_{\mathbf{C}_r} \, \mathbf{x}_r(t).
		\end{aligned}
	\end{equation}
\end{enumerate}
\subsubsection*{Approach 2: Reduction with modal truncation}
\noindent Unfortunately, in this case 3 the application of a one-sided reduction with either an input or an output Krylov subspace would yield time-varying projection matrices $\mathbf{V}(t)$ and $\mathbf{W}(t)$ according to \eqref{equ:V(t)_B(t)} or \eqref{equ:W(t)_C(t)}, respectively. A possible alternative to still obtain time-independent projection matrices and thus get rid of the time derivative $\dot{\mathbf{V}}$ is to use modal truncation as reduction approach. This method only uses the time-independent matrices $\mathbf{A}$ and $\mathbf{E}$ for computing dominant eigenvalues (e.g. with smallest magnitude or smallest real part) and eigenvectors, thus yielding time-independent projection matrices for the reduction. 

\section{Time-Varying Parametric Model Order Reduction} \label{sec:4}
\noindent After having considered linear time-varying systems and presented some straightforward approaches to reduce special cases arising in moving load and sensor problems, in this section we focus on linear parameter-varying systems of the form 

\begin{equation}
\begin{aligned}
\mathbf{E}(\mathbf{p}(t)) \, \dot{\mathbf{x}}(t) &= \mathbf{A}(\mathbf{p}(t)) \, \mathbf{x}(t) + \mathbf{B}(\mathbf{p}(t)) \, \mathbf{u}(t),\\
\mathbf{y}(t) &= \mathbf{C}(\mathbf{p}(t)) \, \mathbf{x}(t).
\end{aligned}
\label{equ:LPV}
\end{equation} 
Such systems also exhibit a time-varying dynamic behaviour, since the system matrices explicitly depend on parameters $\mathbf{p}(t)$ which vary over time. Note that moving load and sensor problems can be represented as LPV systems with only parameter-varying input and/or output matrices, if the load and sensor location are considered as time-dependent parameters of the system model. Due to the time-dependency of the parameters, in the next subsection we derive a projection-based, time-varying parametric model order reduction approach called \emph{p(t)MOR}, to obtain a reduced-order model of a LPV system \cite{Cru_Cruz_MATHMOD_2015,Cru_Cruz_Hirschberg_2015}. Based on that, we then adapt the pMOR approach by matrix interpolation \cite{Cru_Panzer_2010} to the parameter-varying case, whereby new time derivative terms emerge \cite{Cru_Cruz_Hirschberg_2015,Cru_Cruz_MoRePaS_2015}. For the sake of a concise presentation, the time argument $t$ will be omitted in the state, input and output vectors hereafter. 

\subsection{Projective p(t)MOR}
\noindent Similarly as explained in subsection~\ref{subsec:tMOR}, in the case of projection-based time-dependent parametric model order reduction we aim to approximate the state vector $\mathbf{x}$ by $\mathbf{x} \approx \mathbf{V}(\mathbf{p}(t)) \, \mathbf{x}_r$ using a \emph{parameter-varying} projection matrix $\mathbf{V}(\mathbf{p}(t))$. Plugging the corresponding approximation equations for $\mathbf{x}$ and its derivative $\dot{\mathbf{x}}$ in \eqref{equ:LPV}, and applying thereon a properly chosen projection matrix $\mathbf{W}(\mathbf{p}(t))$ that imposes the Petrov-Galerkin condition yields the reduced-order model
\begin{equation}
\begin{aligned}
\mathbf{E}_r(\mathbf{p}(t)) \, \dot{\mathbf{x}}_r &= \left(\mathbf{A}_r(\mathbf{p}(t)) - \mathbf{W}(\mathbf{p}(t))^T \mathbf{E}(\mathbf{p}(t)) \textcolor{tumred}{\dot{\mathbf{V}}(\mathbf{p}(t))} \right) \, \mathbf{x}_r + \mathbf{B}_r(\mathbf{p}(t)) \, \mathbf{u},  \\[0.3em]
\mathbf{y}_r &= \mathbf{C}_r(\mathbf{p}(t)) \, \mathbf{x}_r,
\end{aligned}
\label{equ:p_dot}
\end{equation}
with the time-dependent parametric reduced matrices
\begin{equation}
\resizebox{0.91\linewidth}{!}{
	$\begin{aligned}
	\mathbf{E}_r(\mathbf{p}(t)) &= \mathbf{W}(\mathbf{p}(t))^T \mathbf{E}(\mathbf{p}(t)) \mathbf{V}(\mathbf{p}(t)), & 
	\mathbf{A}_r(\mathbf{p}(t)) &= \mathbf{W}(\mathbf{p}(t))^T \mathbf{A}(\mathbf{p}(t)) \mathbf{V}(\mathbf{p}(t)),\\[0.3em]
	\mathbf{B}_r(\mathbf{p}(t)) &= \mathbf{W}(\mathbf{p}(t))^T \mathbf{B}(\mathbf{p}(t)), &
	\mathbf{C}_r(\mathbf{p}(t)) &= \mathbf{C}(\mathbf{p}(t)) \mathbf{V}(\mathbf{p}(t)).
	\end{aligned}$
}
\end{equation}
The reduced model comprises an additional term depending on the time derivative $\dot{\mathbf{V}}(\mathbf{p}(t))$, which has to be considered during the extension of the matrix interpolation method to the parameter-varying case.

\subsection{p(t)MOR by Matrix Interpolation}
\noindent The local pMOR technique of matrix interpolation can be applied to efficiently obtain a parametric reduced-order model from the interpolation of reduced matrices precomputed at different grid points in the parameter space. Similarly as in the classic method \cite{Cru_Panzer_2010}, the LPV system \eqref{equ:LPV} is first evaluated and individually reduced at certain parameter samples $\mathbf{p}_i, i=1,\ldots,k$ with respective projection matrices $\mathbf{V}_i:=\mathbf{V}(\mathbf{p}_i)$ and $\mathbf{W}_i:=\mathbf{W}(\mathbf{p}_i)$. The reduced state vectors $\mathbf{x}_{r,i}$ of the independently calculated reduced models 
\begin{equation}
\begin{aligned}
\mathbf{E}_{r,i} \, \dot{\mathbf{x}}_{r,i} &= \left(\mathbf{A}_{r,i} - \mathbf{W}_{i}^T \, \mathbf{E}_{i} \, \dot{\mathbf{V}}(\mathbf{p}(t)) \right) \, \mathbf{x}_{r,i} + \mathbf{B}_{r,i} \, \mathbf{u}, \\[0.2em]
\mathbf{y}_{r,i} &= \mathbf{C}_{r,i} \, \mathbf{x}_{r,i}
\end{aligned}
\label{equ:local_models_p(t)}
\end{equation}
generally lie in different subspaces and have, therefore, different physical meanings. For this reason, the direct interpolation of the reduced matrices is not meaningful, and hence the local reduced models have to be transformed into a common set of coordinates first. This is performed applying state transformations of the form
\begin{equation}
\begin{aligned}
\mathbf{x}_{r,i} &= \mathbf{T}_i \, \hat{\mathbf{x}}_{r,i}, \\[0.2em]
\dot{\mathbf{x}}_{r,i} &= \dot{\mathbf{T}}_i \, \hat{\mathbf{x}}_{r,i} + \mathbf{T}_i \, \dot{\hat{\mathbf{x}}}_{r,i},
\end{aligned}
\end{equation}
with regular matrices $\mathbf{T}_i:=\mathbf{T}(\mathbf{p}_i)$, whereby the product rule is required again for the differentiation of $\mathbf{x}_{r,i}$. These state transformations serve to adjust the different \emph{right} local bases $\mathbf{V}_i$ to new bases $\hat{\mathbf{V}}_i=\mathbf{V}_i \, \mathbf{T}_i$. In order to adjust the different \emph{left} local bases $\mathbf{W}_i$ by means of $\hat{\mathbf{W}}_i=\mathbf{W}_i \, \mathbf{M}_i$ as well, the reduced models from \eqref{equ:local_models_p(t)} are subsequently multiplied from the left with regular matrices $\mathbf{M}_i^T$. The resulting reduced and transformed systems are thus given by 
\begin{equation}
\resizebox{0.91\linewidth}{!}{
	$\begin{aligned}
	\overbrace{\mathbf{M}_i^T\mathbf{E}_{r,i}\mathbf{T}_i}^{\hat{\mathbf{E}}_{r,i}} \dot{\hat{\mathbf{x}}}_{r,i}\!&=\!
	\overbrace{\left(\overbrace{\mathbf{M}_i^T\mathbf{A}_{r,i} \mathbf{T}_i}^{\hat{\mathbf{A}}_{r,i}} - \mathbf{M}_i^T\mathbf{W}_i^T \mathbf{E}_i \dot{\mathbf{V}}(\mathbf{p}(t)) \mathbf{T}_i - 
		\mathbf{M}_i^T \mathbf{E}_{r,i} \dot{\mathbf{T}}_i \right)}^{\hat{\mathbf{A}}_{\textrm{new}\,r,i}} \hat{\mathbf{x}}_{r,i}
	+ \overbrace{\mathbf{M}_i^T\mathbf{B}_{r,i}}^{\hat{\mathbf{B}}_{r,i}} \mathbf{u}, \\[0.2em]
	\mathbf{y}_{r,i} \!&=\! \underbrace{\mathbf{C}_{r,i} \mathbf{T}_i}_{\hat{\mathbf{C}}_{r,i}} \hat{\mathbf{x}}_{r,i}.
	\end{aligned}$
}
\label{equ:local_models_transformated_p(t)}
\end{equation}
One possible way to calculate the transformation matrices $\mathbf{T}_i$ and $\mathbf{M}_i$ is based on making the state vectors $\hat{\mathbf{x}}_{r,i}$ compatible with respect to a reference subspace spanned by the columns of the orthogonal matrix $\mathbf{R}$. To this end, the matrices are chosen as $\mathbf{T}_i:=(\mathbf{R}^T \, \mathbf{V}_i)^{-1}$ and $\mathbf{M}_i:=(\mathbf{R}^T \, \mathbf{W}_i)^{-1}$, where the columns of $\mathbf{R}$ correspond to the $r$ most important directions of $\mathbf{V}_{\textrm{all}} = \left[\mathbf{V}_1 \ \ldots \ \mathbf{V}_k\right]$ calculated by a Singular Value Decomposition (SVD) \cite{Cru_Panzer_2010}.

The resulting system matrix $\hat{\mathbf{A}}_{\textrm{new}\,r,i}$ not only comprises the expected reduced matrix $\hat{\mathbf{A}}_{r,i}$, but also consists of two further terms that depend on $\dot{\mathbf{V}}(\mathbf{p}(t))$ and $\dot{\mathbf{T}}_i$, respectively. The calculation of these time derivatives that are required for the computation of the reduced-order model will be discussed in the next two sections.   

After the transformation of the local models and the computation of the new emerging time derivatives, a parameter-varying reduced-order model for a new parameter value $\mathbf{p}(t)$ is obtained in the online phase by a weighted interpolation between the reduced matrices from \eqref{equ:local_models_transformated_p(t)} according to
\begin{equation}
\begin{aligned}
\tilde{\mathbf{E}}_{r}(\mathbf{p}(t)) &= \sum\nolimits_{i=1}^k \omega_{i}(\mathbf{p}(t)) \hat{\mathbf{E}}_{r,i}, & \tilde{\mathbf{A}}_{\textrm{new}\,r}(\mathbf{p}(t)) &= \sum\nolimits_{i=1}^k \omega_{i}(\mathbf{p}(t)) \hat{\mathbf{A}}_{\textrm{new}\,r,i},  \\[0.2em]
\tilde{\mathbf{B}}_{r}(\mathbf{p}(t)) &= \sum\nolimits_{i=1}^k \omega_{i}(\mathbf{p}(t)) \hat{\mathbf{B}}_{r,i}, & \tilde{\mathbf{C}}_{r}(\mathbf{p}(t)) &= \sum\nolimits_{i=1}^k \omega_{i}(\mathbf{p}(t)) \hat{\mathbf{C}}_{r,i},
\end{aligned}
\end{equation}
where $\sum\nolimits_{i=1}^k \omega_{i}(\mathbf{p}(t))=1$. For simplicity, here we use piecewise linear interpolation of the reduced matrices. Higher order interpolation schemes could also be applied.

\subsubsection{Time derivative of $\mathbf{V}$}
\begin{figure}[h]
	\centering
	\tikzsetnextfilename{derivative_V}
	\begin{tikzpicture}

\draw [tumblue3, fill]  (-3, -0.05) rectangle (-1, 0.05);
\draw [tumblue2, fill]  (-1, -0.05) rectangle (1, 0.05); 
\draw [tumblue1, fill]  (1, -0.05) rectangle (3, 0.05);

\draw (-3, 0) node[below=0.6em] {$\underline{\mathbf{p}}_{t-1}$}  +(-0.2, 0.2)  -- +(0.2, -0.2);
\draw (-3, 0) node[above=0.6em] {$\underline{\mathbf{V}}_{t-1}$}  +(-0.2, -0.2)  -- +(0.2, 0.2);

\draw (-1, 0) node[below=0.6em] {$\overline{\mathbf{p}}_{t-1}$}  +(-0.2, 0.2)  -- +(0.2, -0.2);
\draw (-1, 0) node[above=0.6em] {$\overline{\mathbf{V}}_{t-1}$}  +(-0.2, -0.2)  -- +(0.2, 0.2);

\draw (1, 0) node[below=0.6em] {$\underline{\mathbf{p}}_{t}$}  +(-0.2, 0.2)  -- +(0.2, -0.2);
\draw (1, 0) node[above=0.6em] {$\underline{\mathbf{V}}_{t}$}  +(-0.2, -0.2)  -- +(0.2, 0.2);

\draw (3, 0) node[below=0.6em] {$\overline{\mathbf{p}}_{t}$}   +(-0.2, 0.2)  -- +(0.2, -0.2);
\draw (3, 0) node[above=0.6em] {$\overline{\mathbf{V}}_{t}$}  +(-0.2, -0.2)  -- +(0.2, 0.2);

\draw[-latex'] (-4, 0) -- (4, 0) node[below] {$\mathbf{p}$};

\draw[tumblue4] (-2, 0) node[below=0.1em] {$\mathbf{p}_{t-1}$}; 
\draw[tumblue4, fill = white] (-2, -0.8) circle (0.3) node[] {$\frac{\partial \mathbf{V}}{\partial \mathbf{p}}$};
\draw[tumblue4, thick, latex'-] (-2,0) +(0, 0) -- +(0, 2) node[above] {$\mathbf{V}_{t-1}$}; 

\draw[tumblue3, fill = white] (0, -0.8) circle (0.3) node[] {$\frac{\partial \mathbf{V}}{\partial \mathbf{p}}$};

\draw[tumblue2] (2, 0) node[below=0.1em] {$\mathbf{p}_t$};
\draw[tumblue2, fill = white] (2, -0.8) circle (0.3) node[] {$\frac{\partial \mathbf{V}}{\partial \mathbf{p}}$};
\draw[tumblue2, thick, latex'-] (2,0) +(0, 0) -- +(0, 2) node[above] {$\mathbf{V}_{t}$};

\draw[dashed,->] (-2, 1.5) -- (2, 1.5);

\end{tikzpicture}
	\caption{Graphical representation of the calculation of the time derivative $\dot{\mathbf{V}}(\mathbf{p}(t))$ using the local bases $\mathbf{V}_i$ computed at the parameter sample points $\mathbf{p}_i$}
	\label{fig:derivative_V}
\end{figure}
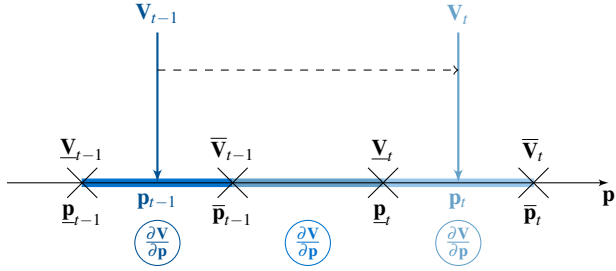

\noindent The time derivative of the projection matrix $\mathbf{V}(\mathbf{p}(t))$ can be numerically calculated using a finite difference approximation. Applying the chain rule first and employing a finite difference method thereon, the time derivative is given by:  
\begin{equation}
\begin{aligned}
\dot{\mathbf{V}}(\mathbf{p}(t)) = \frac{\partial \mathbf{V}}{\partial \mathbf{p}} \, \dot{\mathbf{p}} = \frac{\overline{\mathbf{V}}_t - \underline{\mathbf{V}}_{t}}{\overline{\mathbf{p}}_t - \underline{\mathbf{p}}_{t}} \, \frac{\mathbf{p}_t - \mathbf{p}_{t-1}}{\Delta t}.
\end{aligned}
\end{equation}
$\overline{\mathbf{p}}_t$ and $\underline{\mathbf{p}}_t$ denote the upper and lower limit of the interval $[\underline{\mathbf{p}}_t, \overline{\mathbf{p}}_t]$, in which the parameter vector  $\mathbf{p}_t$ is located at time instant $t$. The local bases at these parameter sample points are given by $\overline{\mathbf{V}}_t$ and $\underline{\mathbf{V}}_t$, respectively. The partial derivatives $\frac{\partial \mathbf{V}}{\partial \mathbf{p}}$ for each pair of parameter sample points are calculated in the offline phase of the matrix interpolation approach. In the online phase, the current time derivative $\dot{\mathbf{V}}(\mathbf{p}(t))$ is then computed by multiplying the partial derivative of the corresponding parameter interval at time instant $t$ with $\dot{\mathbf{p}}$, which represents the current velocity of the moving load. Fig.~\ref{fig:derivative_V} illustrates the aforementioned intervals and the efficient numerical calculation of the time derivative $\dot{\mathbf{V}}(\mathbf{p}(t))$ by a finite difference approximation using only precomputed local bases.
 
\subsubsection{Time derivative of $\mathbf{T}$}
\noindent As explained before, in this paper the transformation matrices $\mathbf{T}_i$ are calculated with $\mathbf{T}_i\!=\!(\mathbf{R}^T \mathbf{V}_i)^{-1}\!:=\!\mathbf{K}^{-1}$. For the computation of the time derivative $\dot{\mathbf{T}}_i$ we make use of the following definition \cite[p.~67]{Cru_Golub_Book_2013}:  

\begin{definition} 
	Let the matrix $\mathbf{K}$ be nonsingular. The time derivative of the inverse matrix is then given by $\frac{\mathrm d \mathbf{K}^{-1}}{\mathrm d t} = -\mathbf{K}^{-1}\frac{\mathrm d \mathbf{K}}{\mathrm d t}\mathbf{K}^{-1}$.
\end{definition}
\noindent This leads to:
\begin{equation}
\begin{aligned}
\dot{\mathbf{T}}_i = \frac{\mathrm d \mathbf{K}^{-1}}{\mathrm d t} = -(\mathbf{R}^T \mathbf{V}_i)^{-1} \mathbf{R}^T \dot{\mathbf{V}}(\mathbf{p}(t)) (\mathbf{R}^T \mathbf{V}_i)^{-1} 
= -\mathbf{T}_i \mathbf{R}^T \dot{\mathbf{V}}(\mathbf{p}(t)) \mathbf{T}_i.
\label{equ:abl_T}
\end{aligned}
\end{equation}

\subsection{p(t)MOR by Matrix Interpolation for particular cases}
\noindent For the reduction of general linear parameter-varying systems the application of time-dependent parametric projection matrices undoubtedly provides an accurate consideration of the arising time variability. Their usage, however, involves some difficulties, like the calculation of the additional derivatives and their incorporation in the numerical simulation of the reduced-order model. Particular LPV systems with only parameter-varying input and/or output matrices, arising e.g. in moving load and sensor problems,  can efficiently be reduced using the matrix interpolation approach combined with the usage of \emph{parameter-independent} projection matrices. In the following, this technique is briefly explained for some special cases:

\subsubsection*{Moving Loads} 
\noindent The application of parameter-varying projection matrices~$\mathbf{V}(\mathbf{p}(t))$ and $\mathbf{W}(\mathbf{p}(t))$ for the individual reduction of the local systems within matrix interpolation results in a reduced model, where \emph{all} reduced matrices vary with the time-dependent parameter, although the original LPV system only contains variations in the input matrix. In order to get rid of the emerging derivatives and \emph{only} have to interpolate the input matrix in the online phase of matrix interpolation, one-sided reductions with an output Krylov subspace $\mathcal{W} = \textrm{span}(\mathbf{W})$ should be employed. 

\subsubsection*{Moving Sensors} 
\noindent In a similar manner, for the case of a LPV system with only parameter-varying output matrix $\mathbf{C}(\mathbf{p}(t))$ one-sided projections with a single input Krylov subspace $\mathcal{V}=\textrm{span}(\mathbf{V})$ computed with the input matrix should be performed for the reduction of the sampled models during matrix interpolation. In this way, we obtain parameter-independent projection matrices $\mathbf{V}=\mathbf{W}$ and only have to interpolate the output matrix, thus reducing the computational effort in the online phase.     

\subsubsection*{Moving Loads and Sensors} 
\noindent For the combined moving load and sensor example the application of one-sided projections with either input or output Krylov subspaces is not helpful, since both the input and output matrices are parameter-varying in this case. Therefore, parameter-independent projection matrices can only be calculated using modal truncation. By doing so, the reduced-order model only contains parameter variations in the input and output reduced matrices like in the original system.

\section{Numerical Examples} \label{sec:5}
\noindent In this section, we present some numerical results for systems with moving loads.

\subsection{Timoshenko beam}
\noindent The presented reduction approaches are first applied to the finite element model of a simply supported Timoshenko beam of lenght $L$ subjected to a moving load.

\begin{figure}[h!]
	\tikzsetnextfilename{beam_moving_load}
	\centering
	\input{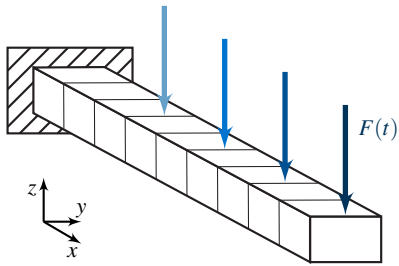}
	\caption{A simply supported Timoshenko beam subjected to a moving force $F(t)$}
	\label{fig:beam_moving_load}
\end{figure}
\noindent Since the moving force $F(t)$ is applied in the negative $z$-direction and we are only interested in the vertical displacement of the beam, the model described in \cite{Cru_Panzer_Beam_2009} is adapted from a 3D to a 1D finite element model. Furthermore, both the moving load and/or sensor case are incorporated into the model, yielding time-dependent input and/or output matrices. The resulting single-input single-output second-order system is reformulated into a LTV first-order model of the form
\begin{equation}
\begin{aligned}
\overbrace{\begin{bmatrix}\mathbf{F} & \mathbf{0}\\\mathbf{0}&\mathbf{M}\end{bmatrix}}^{\mathbf{E}}
\overbrace{\begin{bmatrix}\dot{\mathbf{z}} \\ \ddot{\mathbf{z}}\end{bmatrix}}^{\dot{\mathbf{x}}}(t) &=
\overbrace{\begin{bmatrix}\mathbf{0}&\mathbf{F}\\\mathbf{-K}&\mathbf{-D}\end{bmatrix}}^{\mathbf{A}}
\overbrace{\begin{bmatrix}\mathbf{z}\\\dot{\mathbf{z}}\end{bmatrix}}^{\mathbf{x}}(t) +
\overbrace{\begin{bmatrix}\mathbf{0}\\\hat{\mathbf{b}}(t)\end{bmatrix}}^{\mathbf{b}(t)}F(t), \\[0.2em]
y(t) &= \underbrace{\begin{bmatrix}\hat{\mathbf{c}}(t)^T & \mathbf{0}^T\end{bmatrix}}_{\mathbf{c}(t)^T} \begin{bmatrix}\mathbf{z}\\\dot{\mathbf{z}}\end{bmatrix}(t),
\end{aligned}
\label{equ:LPV_beam}
\end{equation}
where the arbitrary nonsingular matrix $\mathbf{F} \in \mathbb{R}^{2N \times 2N}$ is chosen in our case to $\mathbf{F}=\mathbf{K}$ for the aim of stability preservation using a one-sided reduction (cf. \cite{Cru_Panzer_Beam_2009,Cru_Cruz_Hirschberg_2015}). The dimension of the original model is then given by $n = 2 \cdot 2N$ with $N$ finite elements.

\subsubsection*{Moving load case} 
\begin{figure}[tp!]
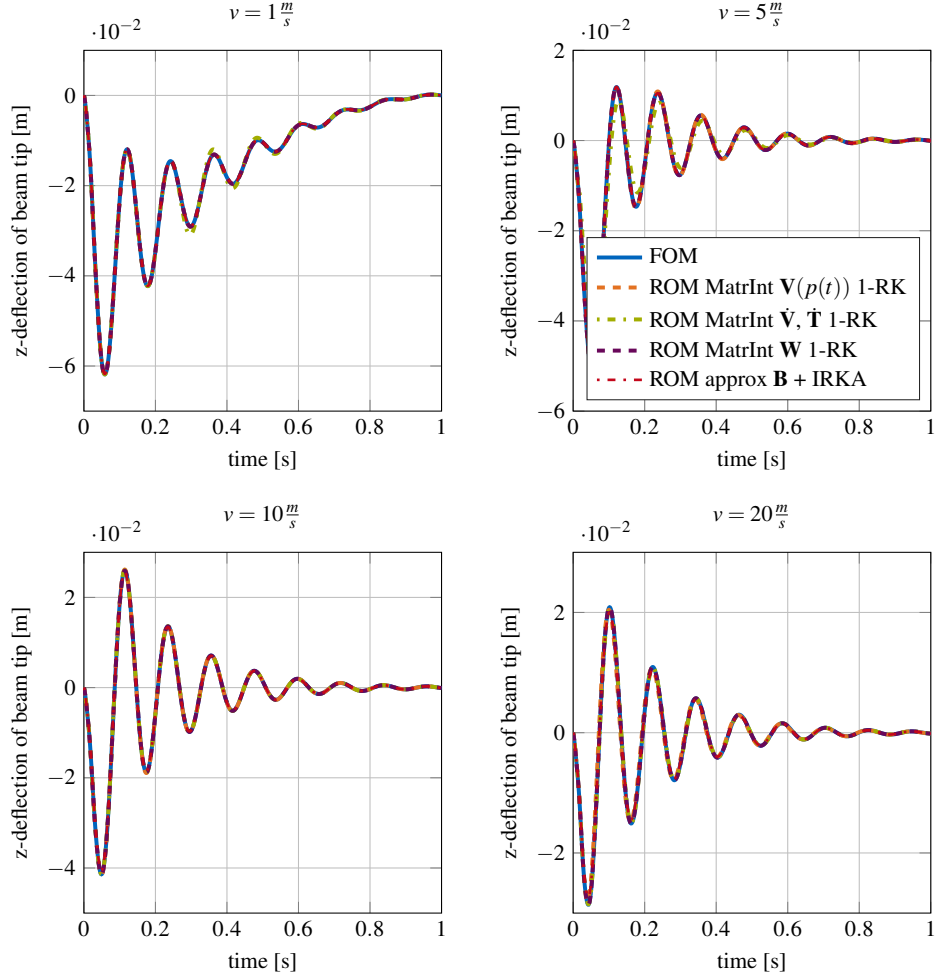

	\centering
	\subfloat{\centering
		\tikzsetnextfilename{movingLoad_L1v1r10} 
		\setlength\myheight{4.8cm}
		\setlength\mywidth{5cm}
		\input{./images/movingLoad_L1v1r10.tikz}		
		\label{fig:movingLoad_L1v1r10}}
	\subfloat{\centering
		\tikzsetnextfilename{movingLoad_L1v5r10} 
		\setlength\myheight{4.8cm}
		\setlength\mywidth{5cm}
		\input{./images/movingLoad_L1v5r10.tikz}		
		\label{fig:movingLoad_L1v5r10}}\\[-0.3em]
	\subfloat{\centering
		\tikzsetnextfilename{movingLoad_L1v10r10} 
		\setlength\myheight{4.8cm}
		\setlength\mywidth{5cm}
		\input{./images/movingLoad_L1v10r10.tikz}		
		\label{fig:movingLoad_L1v10r10}}
	\subfloat{\centering
		\tikzsetnextfilename{movingLoad_L1v20r10} 
		\setlength\myheight{4.8cm}
		\setlength\mywidth{5cm}
		\input{./images/movingLoad_L1v20r10.tikz}		
		\label{fig:movingLoad_L1v20r10}}
	\caption{Simulation results for the Timoshenko beam with moving load for different reduction methods and velocities. Original dimension $n = 2 \cdot 2 \cdot 451 = 1804$, reduced dimension $r=10$. Krylov-based reductions performed with expansion points $s_0=0$}
	\label{fig:movingLoadRed_r10}
\end{figure}

\noindent We first consider the reduction of a beam of length $L=\unit[1]{m}$ subjected to a point force moving from the tip to the supporting with a constant velocity $v$ and an amplitude of $F(t)=\unit[20]{N}$. For the numerical simulation we use an implicit Euler scheme with a step size of $dt = \unit[0.001]{s}$. In Fig. \ref{fig:movingLoadRed_r10} the simulation results for the different proposed reduction methods are presented. We first apply the standard matrix interpolation (MatrInt) approach using $k=76$ equidistantly distributed local models with corresponding current input, which are individually reduced applying one-sided projections with input Krylov subspaces ($\mathbf{V}(p(t))$) for $r=10$. The consideration of the theorically emerging derivatives $\dot{\mathbf{V}}$ and $\dot{\mathbf{T}}$ according to \eqref{equ:local_models_transformated_p(t)} in the matrix interpolation scheme only yields better results than the standard MatrInt method for large velocities of the moving load. In any case, the application of a single \emph{time-independent} output Krylov subspace ($\mathbf{W}$) during MatrInt and the two-step method ($\tilde{m}\!=\!76$) combined with MIMO-IRKA yields the best results (see Table \ref{tab:1}).

\begin{table}
	\caption{Absolute $\mathcal{L}_2$ output error norms $\|y-y_r\|_{\mathcal{L}_2}$}\label{tab:1}
	\begin{tabular}{p{6em}p{7.5em}p{7em}p{7em}p{8em}}
		\hline\noalign{\smallskip}
		& MatrInt $\mathbf{V}(p(t))$ & MatrInt $\dot{\mathbf{V}}$, $\dot{\mathbf{T}}$ & MatrInt $\mathbf{W}$ & approx $\mathbf{B}$ + IRKA \\
		\noalign{\smallskip}\hline\noalign{\smallskip}
		$v=1 \frac{m}{s}$ & $4.8e^{-4}$  & $1.5e^{-2}$ & $3.0e^{-5}$ & $1.0e^{-4}$\\[0.3em]
		$v=5 \frac{m}{s}$ & $3.0e^{-3}$ & $6.8e^{-2}$ & $1.8e^{-5}$ & $2.0e^{-4}$\\[0.3em]
		$v=10 \frac{m}{s}$ & $2.6e^{-3}$  & \textcolor{tumbluegreen}{$1.0e^{-3}$} & $1.7e^{-5}$ & $5.8e^{-5}$\\[0.3em]
		$v=20 \frac{m}{s}$ & $2.2e^{-2}$  & \textcolor{tumbluegreen}{$5.4e^{-3}$} & $1.2e^{-5}$ & $3.9e^{-5}$\\
		\noalign{\smallskip}\hline\noalign{\smallskip}
	\end{tabular}
\end{table}
\subsubsection*{Moving load and sensor case} 
\noindent Now we consider a larger beam of length $L=\unit[50]{m}$ with both moving load and sensor. The observation of the z-deflection of the beam coincides at any time with the position of the moving load, meaning that $\mathbf{c}(p(t))^T = \mathbf{b}(p(t))$. First we apply the matrix interpolation approach and use modal truncation for the individual reduction of the $k=201$ sampled models constructed with the input and output vectors corresponding to each parameter sample point. Since modal truncation only considers the matrices $\mathbf{A}$ and $\mathbf{E}$ for the reduction and these matrices do not vary over time, we only have to compute \emph{one single pair} of time-independent projection matrices $\mathbf{V}$ and $\mathbf{W}$ in the offline phase. During the online phase, only the parameter-varying input and output vectors have to be interpolated in order to obtain a reduced-order model for each current position of the load/sensor. 

\begin{figure}[htp!]
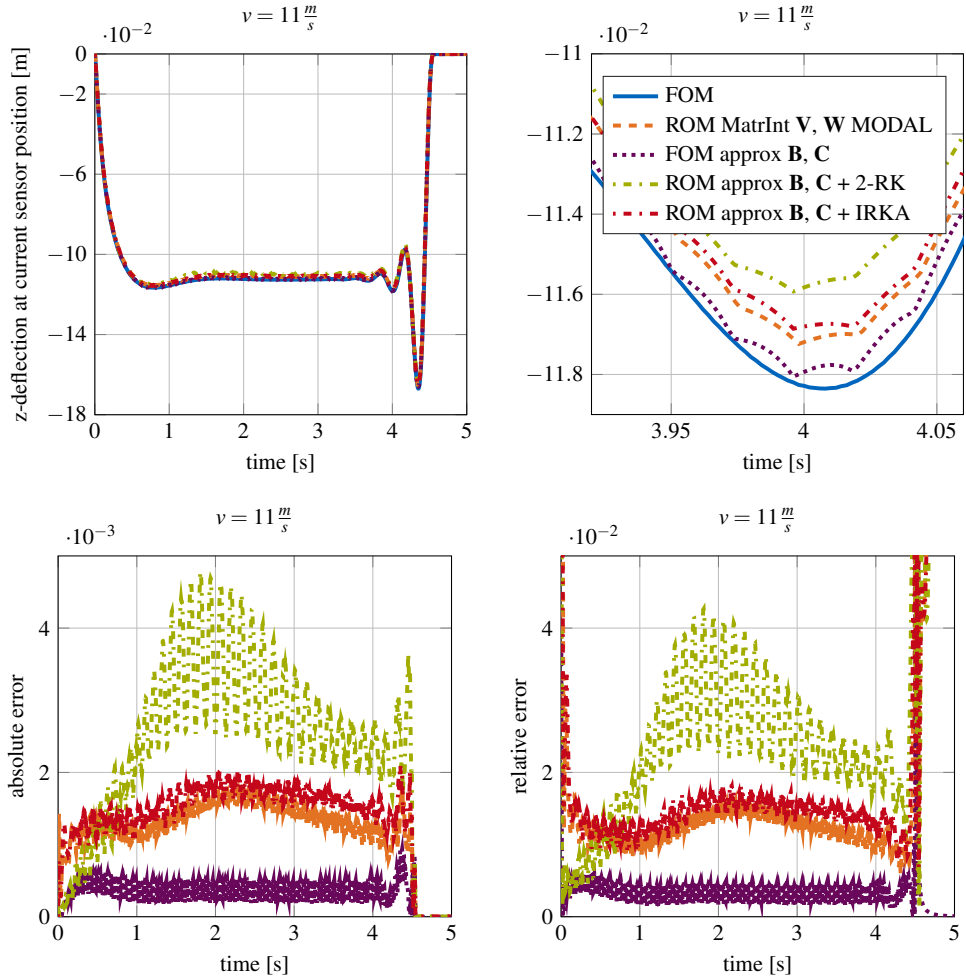

	\centering
	\subfloat{\centering
		\tikzsetnextfilename{movingLoadSensor_L50v11r80_y} 
		\setlength\myheight{4.8cm}
		\setlength\mywidth{5.2cm}
		\input{./images/movingLoadSensor_L50v11r80_y.tikz}		
		\label{fig:movingLoadSensor_L50v11r80_y}}
	\subfloat{\centering
		\tikzsetnextfilename{movingLoadSensor_L50v11r80_yZoom} 
		\setlength\myheight{4.8cm}
		\setlength\mywidth{5.2cm}
		\input{./images/movingLoadSensor_L50v11r80_yZoom.tikz}		
		\label{fig:movingLoadSensor_L50v11r80_yZoom}}\\[-0.3em]
	\subfloat{\centering
		\tikzsetnextfilename{movingLoadSensor_L50v11r80_absError} 
		\setlength\myheight{4.8cm}
		\setlength\mywidth{5.5cm}
		\input{./images/movingLoadSensor_L50v11r80_absError.tikz}		
		\label{fig:movingLoadSensor_L50v11r80_absError}}
	\subfloat{\centering
		\tikzsetnextfilename{movingLoadSensor_L50v11r80_relError} 
		\setlength\myheight{4.8cm}
		\setlength\mywidth{5.5cm}
		\input{./images/movingLoadSensor_L50v11r80_relError.tikz}		
		\label{fig:movingLoadSensor_L50v11r80_relError}}
	\caption{Simulation results for the Timoshenko beam with moving load and sensor for different reduction methods. Original dimension $n = 2 \cdot 2 \cdot 1001 = 4004$, reduced dimension $r=80$. Krylov-based reductions performed with expansion points $s_0=0$}
	\label{fig:movingLoadSensor_L50v11r80}
\end{figure}

Next, we further apply the aforementioned two-step method for the reduction. To this end, the time-varying input and output vectors are first approximated by low-rank matrices $\mathbf{B}$ and $\mathbf{C}$ on a coarse finite element grid. To ensure a proper comparability with MatrInt, we choose the same $\tilde{m}=201$ nodes where local models were constructed before. The herewith obtained approximated output $y(t)$ and approximation errors are depiced in Fig. \ref{fig:movingLoadSensor_L50v11r80}. One can see that the number of chosen columns $\tilde{m}$ is sufficiently large, since the approximation error is adequately small. After that, we both apply two-sided MIMO rational Krylov (2-RK) and MIMO-IRKA for the reduction of the resulting LTI system. Fig.~\ref{fig:movingLoadSensor_L50v11r80} shows the simulated output for the different explained reduction methods as well as the corresponding absolute and relative $\mathcal{L}_2$ errors. Although all results show a similar behaviour, the matrix interpolation approach combined with modal truncation together with the two-step method by IRKA lead to the smallest errors. Simulations were also conducted with the extended p(t)MOR approach by matrix interpolation considering the time derivatives like in \ref{equ:local_models_transformated_p(t)}. Unfortunately, these additional terms make the pencils $(\hat{\mathbf{A}}_{\textrm{new}\,r,i}, \hat{\mathbf{E}}_{r,i})$ often unstable, yielding unstable interpolated systems and results.   

\subsection{Beam with moving heat source} \label{subsec:1D_heat}
\begin{figure}[tp!]
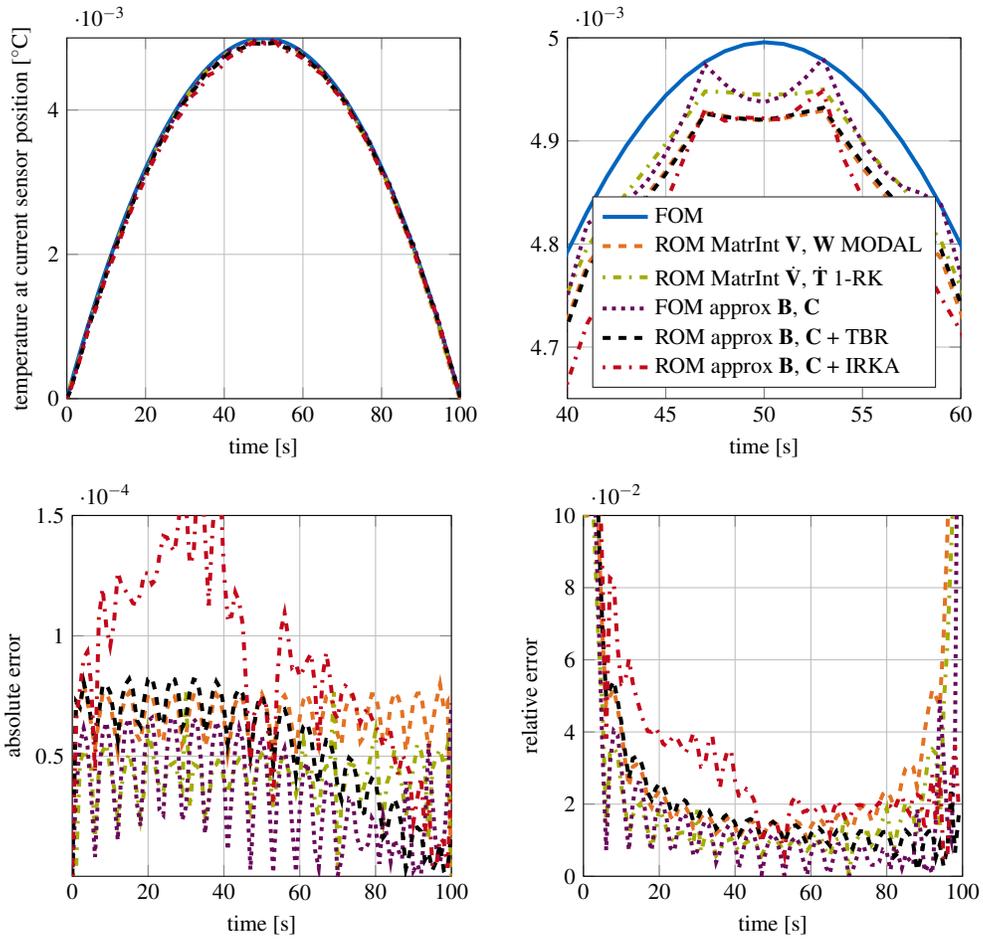
 
	\centering
	\subfloat{\centering
		\tikzsetnextfilename{1DHeat_movingLoadSensor_L10r40_y} 
		\setlength\myheight{4.8cm}
		\setlength\mywidth{5.5cm}
		\input{./images/1DHeat_movingLoadSensor_L10r40_y.tikz}		
		\label{fig:1DHeat_movingLoadSensor_L10r40_y}}
	\subfloat{\centering
		\tikzsetnextfilename{1DHeat_movingLoadSensor_L10r40_yZoom} 
		\setlength\myheight{4.8cm}
		\setlength\mywidth{5.5cm}
		\input{./images/1DHeat_movingLoadSensor_L10r40_yZoom.tikz}		
		\label{fig:1DHeat_movingLoadSensor_L10r40_yZoom}}\\[-0.3em]
	\subfloat{\centering
		\tikzsetnextfilename{1DHeat_movingLoadSensor_L10r40_absError} 
		\setlength\myheight{4.8cm}
		\setlength\mywidth{5.3cm}
		\input{./images/1DHeat_movingLoadSensor_L10r40_absError.tikz}		
		\label{fig:1DHeat_movingLoadSensor_L10r40_absError}}
	\subfloat{\centering
		\tikzsetnextfilename{1DHeat_movingLoadSensor_L10r40_relError} 
		\setlength\myheight{4.8cm}
		\setlength\mywidth{5.3cm}
		\input{./images/1DHeat_movingLoadSensor_L10r40_relError.tikz}		
		\label{fig:1DHeat_movingLoadSensor_L10r40_relError}}
	\caption{Simulation results for the 1D beam with moving heat source for different reduction methods. Original dimension $n=2500$, reduced dimension $r=40$. Krylov-based reductions performed with expansion points $s_0=0$}
	\label{fig:1DHeat_movingLoadSensor_L10r40}
\end{figure}

\noindent We now apply the presented techniques on a second example \cite{Cru_Lang_MCMDS_2016}, which describes the heat transfer along a beam of length $L$ with a moving heat source. The temperature is observed at the same position as the heat source, thus $\mathbf{c}(t)^T = \mathbf{b}(t)$. In our case, we consider a system dimension of $n=2500$, apply an input heating source of $u(t)=\unit[50]{^{\circ} C}$ and use an implicit Euler scheme with a step size of $dt = \unit[1]{s}$ for the time integration. Fig. \ref{fig:1DHeat_movingLoadSensor_L10r40} shows the simulation results, and the absolute and relative errors for the different employed reduction methods. One interesting observation is that in this case the application of the extended MatrInt approach with the consideration of the time derivatives yields a slightly better approximation than the classic MatrInt combined with modal truncation ($k=84$). In general, this fact could also be observed for the previous and some other numerical experiments with higher velocities, as long as the overall interpolated systems were stable. This slightly better approximation can be explained through the more accurate consideration of the arising time variability using time-dependent projection matrices, as opposed to modal truncation which does not consider the moving interactions. The approximation of the time-dependent input and output vectors by low-rank matrices using $\tilde{m}=84$ nodes, and the subsequent application of balanced truncation (TBR) or MIMO-IRKA for the reduction shows a similar behaviour. Although the extended MatrInt shows in this case the best results, it is difficult to clearly identify a superior method, since all presented approaches are suitable for the reduction of systems with moving loads.  

\section{Conclusions} \label{sec:6}
\noindent In this paper, we have presented several time-varying model reduction techniques for systems with moving loads. Such particular, but still frequent problems lead to high-dimensional systems with time-varying input and/or output matrices. For their reduction, time-dependent projection matrices can be applied, thus offering an accurate consideration of the time variation, but leading also to an additional derivative in the reduced model which have to be taken into account. Since moving load problems represent particular LTV systems, we have presented straightforward reduction approaches for some special cases, where time-independent projection matrices are calculated and therefore the emerging time derivative is avoided. Systems with moving loads can also be modeled as special LPV systems, where the input and/or output matrices depend on a time-varying parameter describing the position of the load. In this context we have derived a projection-based, time-varying parametric model reduction approach and extended the matrix interpolation scheme to the parameter-varying case. With the appropriate combination of this method with the application of parameter-independent projection matrices, special LPV systems can be efficiently reduced avoiding the time derivatives. The proposed methods have been tested on two different beam models for both the moving load and/or sensor cases. All techniques have provided similar satisfactory results, showing that all methods are suitable for the reduction of systems with moving loads. In particular, the presented straightforward approaches using time-independent projection matrices are very simple, but may be absolutely sufficient for certain problems. They provide a basis for comparison with more complex techniques that consider the time variability using time-dependent projection matrices. These advanced techniques should be investigated more deeply in the future, especially concerning \emph{general} LTV systems, the increased computational effort due to the time-dependent projection matrices and derivatives, fast-varying load variations and stability preservation in the reduced-order model.

\ack 
\noindent The authors would like to thank N. Lang, T. Stykel and J. Saak for kindly proving us the 1D heat transfer model used for the numerical example in section \ref{subsec:1D_heat}. Furthermore, we thank the former and current members of our model reduction lab for the fruitful discussions.

\bibliography{ref}

\begin{thebibliography}{10}

\bibitem{Cru_Antoulas_Book}
A.~C. Antoulas.
\newblock {\em Approximation of Large-Scale Dynamical Systems}.
\newblock SIAM, Philadelphia, PA, 2005.

\bibitem{Cru_Baumann_Comparison_2016}
M.~Baumann, A.~Vasilyev, T.~Stykel, and P.~Eberhard.
\newblock Comparison of two model order reduction methods for elastic multibody
  systems with moving loads.
\newblock In {\em Proceedings of the Institution of Mechanical Engineers, Part
  K: Journal of Multi-body Dynamics}, 2016.

\bibitem{Cru_Baur_2014_Survey}
U.~Baur, P.~Benner, and L.~Feng.
\newblock Model order reduction for linear and nonlinear systems: a
  system-theoretic perspective.
\newblock {\em Archives of Computational Methods in Engineering},
  21(4):331--358, 2014.

\bibitem{Cru_Benner_2015_Survey}
P.~Benner, S.~Gugercin, and K.~Willcox.
\newblock A survey of projection-based model reduction methods for parametric
  dynamical systems.
\newblock {\em SIAM Review}, 57(4):483--531, 2015.

\bibitem{Cru_Cruz_MATHMOD_2015}
M.~Cruz~Varona, M.~Geuss, and B.~Lohmann.
\newblock p(t){MOR}: Time-varying parametric model order reduction and
  applications for moving loads.
\newblock In {\em Proceedings of the 8th Vienna Conference on Mathematical
  Modelling (MATHMOD)}, volume~48, pages 677--678. Elsevier, 2015.

\bibitem{Cru_Cruz_Hirschberg_2015}
M.~Cruz~Varona, M.~Geuss, and B.~Lohmann.
\newblock Zeitvariante parametrische {M}odellordnungsreduktion am {B}eispiel
  von {S}ystemen mit wandernder {L}ast.
\newblock In {\em Methoden und Anwendungen der Regelungstechnik --
  Erlangen-M{\"u}nchener Workshops 2013 und 2014}, pages 57--70. B. Lohmann und
  G. Roppenecker (Hrsg.), Shaker-Verlag, 2015.

\bibitem{Cru_Cruz_MoRePaS_2015}
M.~Cruz~Varona and B.~Lohmann.
\newblock Time-varying parametric model order reduction by matrix
  interpolation.
\newblock In {\em Proceedings of the MoRePaS 2015--Model Reduction of
  Parametrized Systems III}, 2015.

\bibitem{Cru_Fischer_2014_Application}
M.~Fischer and P.~Eberhard.
\newblock Application of parametric model reduction with matrix interpolation
  for simulation of moving loads in elastic multibody systems.
\newblock {\em Advances in Computational Mathematics}, pages 1--24, 2014.

\bibitem{Cru_Fischer_ECCOMAS_2015}
M.~Fischer and P.~Eberhard.
\newblock Interpolation-based parametric model order reduction for material
  removal in elastic multibody systems.
\newblock In {\em Proceedings ECCOMAS Thematic Conference on Multibody
  Dynamics}, 2015.

\bibitem{Cru_Golub_Book_2013}
G.~H. Golub and C.~F. Van~Loan.
\newblock {\em Matrix Computations}.
\newblock Johns Hopkins University Press, Baltimore, 4th edition, 2013.

\bibitem{Cru_Lang_2014}
N.~Lang, J.~Saak, and P.~Benner.
\newblock Model order reduction for systems with moving loads.
\newblock {\em at-Automatisierungstechnik}, 62(7):512--522, June 2014.

\bibitem{Cru_Lang_MCMDS_2016}
N.~Lang, J.~Saak, and T.~Stykel.
\newblock Balanced truncation model reduction for linear time-varying systems.
\newblock {\em MCMDS}, to appear.

\bibitem{Cru_Panzer_Beam_2009}
H.~Panzer, J.~Hubele, R.~Eid, and B.~Lohmann.
\newblock Generating a parametric finite element model of a {3D} cantilever
  {Timoshenko} beam using {MATLAB}.
\newblock {TRAC}, Lehrstuhl {f\"ur} Regelungstechnik, {Technische Universit\"at
  M\"unchen}, 09 2009.

\bibitem{Cru_Panzer_2010}
H.~Panzer, J.~Mohring, R.~Eid, and B.~Lohmann.
\newblock Parametric model order reduction by matrix interpolation.
\newblock {\em at--Automatisierungstechnik}, 58(8):475--484, 8 2010.

\bibitem{Cru_Sandberg_2004_LTV-BT}
H.~Sandberg and A.~Rantzer.
\newblock Balanced truncation of linear time-varying systems.
\newblock {\em IEEE Transactions on Automatic Control}, 49(2):217--229, 2004.

\bibitem{Cru_Shokoohi_1983_LTV-BT}
S.~Shokoohi, L.~M. Silverman, and P.~M. van Dooren.
\newblock Linear time-variable systems: balancing and model reduction.
\newblock {\em IEEE Transactions on Automatic Control}, 28(8):810--822, 1983.

\bibitem{Cru_Stykel_CAM_2016}
T.~Stykel and A.~Vasilyev.
\newblock A two-step model reduction approach for mechanical systems with
  moving loads.
\newblock {\em Journal of Computational and Applied Mathematics}, 297:85--97,
  2016.

\bibitem{Cru_Tamarozzi_2014}
T.~Tamarozzi, G.~Heirman, and W.~Desmet.
\newblock An on-line time dependent parametric model order reduction scheme
  with focus on dynamic stress recovery.
\newblock {\em Computer Methods in Applied Mechanics and Engineering},
  268(0):336 -- 358, 2014.

\end{thebibliography}
\bibliographystyle{plainnat}

\end{document}